\documentclass[11pt]{article}

\usepackage{amsmath}
\usepackage{amsfonts}
\usepackage{amssymb}
\usepackage{pst-all}
\input epsf

\newtheorem{theorem}{Theorem}

\newtheorem{lemma}{Lemma}
\newtheorem{proposition}{Proposition}

\newenvironment{proof}[1][Proof]{\noindent\textbf{#1.} }{\ \rule{0.5em}{0.5em}}

\def\I{\infty}

\def\l2{\mathcal{L}_{q,2,v}}

\begin{document}

\title{\textbf{On the Graf's addition theorem for
Hahn Exton q-Bessel function}}
\author{Lazhar Dhaouadi\thanks{
IPEB, 7021 Zarzouna, Bizerte,Tunisia. E-mail
lazhardhaouadi@yahoo.fr}\quad  Ahmed Fitouhi\thanks{ Facult\'e des
sciences de Tunis, 1060 Tunis, Tunisia. E-mail
ahmed.fitouhi@fst.rnu.tn}}
\date{ }
\maketitle

\begin{abstract}
In this paper we study the positivity of the generalized
$q$-translation associated with the $q$-Bessel Hahn Exton function
which is deduced by a new formulation of the Graf's addition formula
related to this function.
\end{abstract}

\section{Introduction and Preliminaries}

It is well known that the generalized translation operator $T$ associated
with the Bessel function of the first kind is positive in the sense if $f>0$%
, then $Tf>0$. This property is easily seen when we write $Tf$ as an
integral representation with a kernel involving the area of some
triangle ([1],[9]) and has several applications in many mathematical
fields such that hypergroup structure, heat equation....\newline In
2002, the appropriated $q$-generalized translation for the
$q$-Bessel Hahn Exton function was founded [3] and the problem of
its positivity asked. In literature we meet many attempts to show
this property in particular case, nerveless they are no definitive
response at to day. In [2] the authors prove that for the
$q$-cosinus , the correspondent $q$-translation operator is positive
if $q\in [0, q_0]$ for some $q_0$. In this work and owing a new
formulation of the Graf's addition theorem [7] we give an
affirmative answer about this theme by a technic involving some
inclusion of sets.\newline To make this work self containing, we
begin by the following preliminaries. Throughout this paper we
consider $0<q<1$ and we adopt the standard conventional notations of
[4]. For complex $a$ We put

\begin{equation*}
(a;q)_0=1,\quad (a;q)_n=\prod_{i=0}^{n-1}(1-aq^{i}), \quad n=1...\infty.
\end{equation*}
Jackson's $q$-integral (see [5]) over the interval $[0,\infty
\lbrack $ is defined by
\begin{equation*}
\int_{0}^{\infty }f(x)d_{q}x=(1-q)\sum_{n=-\infty }^{\infty }q^{n}f(q^{n}).
\end{equation*}
We denote by

\begin{equation*}
\mathbb{R}_{q}^{+}=\{q^{n},\quad n\in \mathbb{Z}\},
\end{equation*}%
and we consider $\mathcal{L}_{q,p,v}$ the space of even functions $f$
defined on $\mathbb{R}_{q}^{+}$ such that
\begin{equation*}
\Vert f\Vert _{q,p,v}=\left[ \int_{0}^{\infty }|f(x)|^{p}x^{2v+1}d_{q}x%
\right] ^{1/p}<\infty .
\end{equation*}%
The $q$-Bessel function of third kind and of order $v$, called also
Hahn-Exton function,is defined by the $q$-hypergeometric function
(see [8])
\begin{equation*}
J_{v}(x,q)=\frac{(q^{v+1};q)_{\infty }}{(q;q)_{\infty }}x^{v}{_{1}\phi _{1}}%
(0,q^{v+1},q,qx^{2}),\quad \Re (v)>-1,
\end{equation*}%
and has a normalized form is given by
\begin{equation*}
j_{v}(x,q)=\frac{(q;q)_{\infty }}{(q^{v+1};q)_{\infty }}x^{-v}J_{v}(x,q)={%
_{1}\phi _{1}}(0,q^{v+1},q,qx^{2})\newline
=\sum_{n=0}^{\infty }(-1)^{n}\frac{q^{\frac{n(n+1)}{2}}}{%
(q;q)_{n}(q^{v+1};q)_{n}}x^{2n}.
\end{equation*}%
It's an entire analytic function in $z$.\newline In([2],[6]) the
q-Bessel Fourier transform $\mathcal{F}_{q,v}$ and its properties
were studied in great detail, it is defined as follow
\begin{equation*}
\mathcal{F}_{q,v}f(x)=c_{q,v}\int_{0}^{\infty
}f(t)j_{v}(xt,q^{2})t^{2v+1}d_{q}t,
\end{equation*}%
where
\begin{equation*}
c_{q,v}=\frac{1}{1-q}\frac{(q^{2v+2},q^{2})_{\infty }}{(q^{2},q^{2})_{\infty
}}.
\end{equation*}%
There is many way to define the $q$-Bessel translation operator [2
],[3]. One of them can be enounced for suitable function $f$ as
follows:
\begin{equation*}
T_{q,x}^{v}f(y)=c_{q,v}\int_{0}^{\infty }\mathcal{F}%
_{q,v}(f)(t)j_{v}(xt,q^{2})j_{v}(yt,q^{2})t^{2v+1}d_{q}t,\quad \forall
x,y\in \mathbb{R}_{q}^{+},\forall f\in \mathcal{L}_{q,1,v}.
\end{equation*}%
Now we say that $T_{q,x}^{v}$ is positive if $T_{q,x}^{v}f\geq 0$ for $f\geq
0$.\newline
Let us putting the domain of positivity of $T_{q,x}^{v}$ by
\begin{equation*}
Q_{v}=\{q\in ]0,1[,\quad T_{q,x}^{v}\quad \text{is positive for all}\quad
x\in \mathbb{R}_{q}^{+}\}.
\end{equation*}%
Trough the result of  [2], $Q_{v}$ is not empty. So for $q\in Q_{v}$
the $q-$convolution product of the two functions $f$ and $g\in \mathcal{L}%
_{q,1,v}$ is defined by
\begin{equation*}
f\ast _{q}g(x)=c_{q,v}\int_{0}^{\infty }T_{q,x}^{v}f(y)g(y)y^{2v+1}d_{q}y.
\end{equation*}%
To close this section we present the following results proved in [2]
which will be used in the remainder.

\begin{proposition}
\begin{equation*}
|j_{v}(q^{n},q^{2})|\leq \frac{(-q^{2};q^{2})_{\I}(-q^{2v+2};q^{2})_{\infty }%
}{(q^{2v+2};q^{2})_{\infty }}\left\{
\begin{array}{c}
1\quad \quad \quad \quad \quad \text{if}\quad n\geq 0 \\
q^{n^{2}-(2v+1)n}\quad \text{if}\quad n<0%
\end{array}%
\right. .
\end{equation*}
\end{proposition}

\begin{theorem}
The operator $\mathcal{F}_{q,v}$ satisfies

\bigskip

1. For all functions $f\in\mathcal{L}_{q,1,v}$,%
\begin{equation*}
\mathcal{F}_{q,v}^2f(x)=f(x),\quad\forall x\in\mathbb{R}_q^+.
\end{equation*}

2. If $f\in \mathcal{L}_{q,1,v}$ and $\mathcal{F}_{q,v}f\in \mathcal{L}%
_{q,1,v}$ then
\begin{equation*}
\Vert \mathcal{F}_{q,v}f\Vert _{q,v,2}=\Vert f\Vert _{q,v,2}.
\end{equation*}
\end{theorem}

\section{The Graf's addition formula}

The Graf's addition formula for Hahn-Exton $q$-Bessel function
proved by H.T.Koelink and F. Swarttouw [7] plays a central role here
. It can be stated as follows.
$$\aligned
&J_{v}(Rq^{1/2(y+z+v)},q)J_{x-v}(q^{1/2z},q)\\
&=\sum_{k\in \mathbb{Z}%
}J_{k}(Rq^{1/2(x+y+k)},q)J_{v+k}(Rq^{1/2(y+k+v)},q)J_{x}(q^{1/2(z-k)},q);
\endaligned$$
and it is valid when $z\in \mathbb{Z}$ and $R,x,y,v\in
\mathbb{C}$ satisfying
\begin{equation*}
|R|^{2}q^{1+\Re (x)+\Re (y)}<1,\quad \Re (x)>-1,\quad R\neq 0.
\end{equation*}%
This formula has originally been derived for $v,x,y\in \mathbb{Z}$,
$R>0$ by the interpretation of the Hahn-Exton $q$-Bessel function as
matrix elements of irreducible unitary representation of the quantum
group of plane motions.\newline If we replace $q$ by $q^{2}$ and $R$
by $q^{r}$ in the previous formula we get:
\begin{align*}
& J_{v}(q^{y+z+v+r},q^{2})J_{x-v}(q^{z},q^{2}) \\
& =\sum_{k\in \mathbb{Z}%
}J_{k}(q^{x+y+k+r},q^{2})J_{v+k}(q^{y+k+v+r},q^{2})J_{x}(q^{z-k},q^{2}),
\end{align*}%
and put
\begin{equation*}
m=y+z+v+r,
\end{equation*}%
so
\begin{equation*}
x+y+k+r=m+k+x-z-v,
\end{equation*}%
\begin{equation*}
y+k+v+r=m+k-z,
\end{equation*}%
and we have
\begin{align*}
& J_{v}(q^{m},q^{2})J_{x-v}(q^{z},q^{2}) \\
& =\sum_{k\in \mathbb{Z}%
}J_{k}(q^{m+k+x-z-v},q^{2})J_{v+k}(q^{m+k-z},q^{2})J_{x}(q^{z-k},q^{2}).
\end{align*}%
This last formula is valid for $z\in \mathbb{Z}$ and $r,x,y,v\in \mathbb{C}$
satisfying
\begin{equation*}
{1+2\Re (r)+\Re (x)+\Re (y)}=1+\Re (r)+\Re (m)-\Re (z)-\Re (v)>0,\quad \Re
(x)>-1.
\end{equation*}%
In the above sum we replace $z-k$ by $k$ we get
\begin{align*}
& J_{v}(q^{m},q^{2})J_{x-v}(q^{z},q^{2}) \\
& =\sum_{k\in \mathbb{Z}%
}J_{z-k}(q^{m+x-v-k},q^{2})J_{v+z-k}(q^{m-k},q^{2})J_{x}(q^{k},q^{2}),
\end{align*}%
The sum in the second member exists for
\begin{equation*}
\forall z\in \mathbb{Z},\quad \forall m,v,x\in \mathbb{C},\quad \Re (x)>-1.
\end{equation*}%
In fact there exist an infinity complex number $r\in \mathbb{C}$ for which
\begin{equation*}
1+\Re (r)+\Re (m)-\Re (z)-\Re (v)>0.
\end{equation*}%
Now using the definition of the normalized $q$-Bessel function
\begin{equation*}
J_{x}(q^{k},q^{2})=\frac{(q^{2x+2},q^{2})_{\infty }}{(q^{2},q^{2})_{\infty }}%
q^{vk}j_{x}(q^{k},q^{2})=(1-q)c_{q,x}q^{xk}j_{x}(q^{k},q^{2}),
\end{equation*}%
we obtain
\begin{align*}
&
q^{mv+z(x-v)}(1-q)^{2}c_{q,v}c_{q,x-v}j_{v}(q^{m},q^{2})j_{x-v}(q^{z},q^{2})
\\
& =(1-q)c_{q,x}\sum_{k\in \mathbb{Z}%
}q^{xk}J_{z-k}(q^{x-v+m-k},q^{2})J_{v+z-k}(q^{m-k},q^{2})j_{x}(q^{k},q^{2}).
\end{align*}%
Let
\begin{equation*}
\lambda =q^{n},\quad n\in \mathbb{Z},
\end{equation*}%
and replace
\begin{equation*}
k\rightarrow k+n,m\rightarrow m+n,z\rightarrow z+n.
\end{equation*}%
This implies
\begin{align*}
& q^{mv+z(x-v)}(1-q)^{2}c_{q,v}c_{q,x-v}j_{v}(q^{m}\lambda
,q^{2})j_{x-v}(q^{z}\lambda ,q^{2}) \\
& =(1-q)c_{q,x}\sum_{k\in \mathbb{Z}%
}q^{xk}J_{z-k}(q^{x-v+m-k},q^{2})J_{v+z-k}(q^{m-k},q^{2})j_{x}(q^{k}\lambda
,q^{2}) \\
& =(1-q)c_{q,x}\sum_{k\in \mathbb{Z}%
}q^{2k(x+1)}q^{-k(x+2)}J_{z-k}(q^{x-v+m-k},q^{2})J_{v+z-k}(q^{m-k},q^{2})j_{x}(q^{k}\lambda ,q^{2}).
\end{align*}%
We put
\begin{equation*}
E_{v,x}(q^{m},q^{z},q^{k})=\frac{1}{(1-q)^{2}c_{q,v}c_{q,x-v}}%
q^{-k(x+2)-mv-z(x-v)}J_{z-k}(q^{x-v+m-k},q^{2})J_{v+z-k}(q^{m-k},q^{2}).
\end{equation*}%
The Proposition 1 shows that the function:
\begin{equation*}
\lambda \mapsto j_{v}(q^{m}\lambda ,q^{2})j_{x-v}(q^{z}\lambda ,q^{2}),\quad
\forall m,z\in \mathbb{Z},
\end{equation*}%
belongs to the space $\mathcal{L}_{q,1,v}$. Theorem 1, part 1 leads to the
following statement:

\begin{proposition}
For $z,m\in\mathbb{Z}$ and $x,v\in\mathbb{C}$ satisfying
\begin{equation*}
\Re(x-v)>-1,\quad\Re(v)>-1
\end{equation*}
we have
\begin{equation*}
j_{v}(q^{m}\lambda,q^{2})j_{x-v}(q^{z}\lambda,q^{2})=c_{q,x}\int_{0}^{\infty
}E_{v,x}(q^{m},q^{z},t)j_{x}(\lambda t,q^{2})t^{2x+1}d_{q}t,
\end{equation*}
and
\begin{equation*}
E_{v,x}(q^{m},q^{z},q^{k})=c_{q,x}\int_{0}^{\infty}j_{v}(q^{m}\lambda
,q^{2})j_{x-v}(q^{z}\lambda,q^{2})j_{x}(q^{k}\lambda,q^{2})%
\lambda^{2x+1}d_{q}\lambda.
\end{equation*}
\end{proposition}

\section{Positivity of the $q$-translation operator}

This section is a direct application of the previous one but before
any things we recall that the $q$-translation operator possess the
$q$-integral representation (see [2])

\begin{proposition}
Let $f\in \mathcal{L}_{q,1,v}$ then
\begin{equation*}
T_{q,x}^{v}f(y)=\int_{0}^{\infty }f(z)D_{v}(x,y,z)z^{2v+1}d_{q}z,
\end{equation*}%
where
\begin{equation*}
D_{v}(x,y,z)=c_{q,v}^{2}\int_{0}^{\infty
}j_{v}(xt,q^{2})j_{v}(yt,q^{2})j_{v}(zt,q^{2})t^{2v+1}d_{q}t.
\end{equation*}
\end{proposition}

When $q$ tends to $1^{-}$ we obtain at least formally the classical
one and the $q$-kernel $D_{v}(x,y,z)$ tends to the classical one
([1],[9]) which involves the area of some triangle. Hence the
positivity of $T_{q,x}^{v}$ is subject of those of $D_{v}(x,y,z)$. A
direct consequence of the above result is the fact that if the
operator $T_{q,x}^{v}$ is positive and $f\in \mathcal{L}_{q,1,v}$
then $T_{q,x}^{v}f\in \mathcal{L}_{q,1,v}$ ( in general it is not
true without this hypothesis). Indeed
$$\aligned
&\int_{0}^{\infty }|T_{q,x}^{v}f(y)|y^{2v+1}d_{q}t \\
&\leq\int_{0}^{\infty }T_{q,x}^{v}|f|(y)y^{2v+1}d_{q}t
=\int_{0}^{\infty }|f(z)|\left[ \int_{0}^{\infty }D_{v}(x,y,z)y^{2v+1}d_{q}y%
\right] z^{2v+1}d_{q}z.
\endaligned$$
Putting
\begin{equation*}
\phi :t\mapsto j_{v}(xt,q^{2})j_{v}(zt,q^{2}),
\end{equation*}%
then we can write
\begin{equation*}
D_{v}(x,y,z)=c_{q,v}\mathcal{F}_{q,v}\phi (y).
\end{equation*}%
This gives by the of the inversion formula in Theorem 1 the important result
as the classical; one
\begin{equation*}
\int_{0}^{\infty }D_{v}(x,y,z)y^{2v+1}d_{q}y=\mathcal{F}_{q,v}^{2}\phi
(0)=\phi (0)=1.
\end{equation*}%
Then we have
\begin{equation*}
\int_{0}^{\infty }|T_{q,x}^{v}f(y)|y^{2v+1}d_{q}t\leq \Vert f\Vert _{q,v,1},
\end{equation*}%
and we obtain
\begin{equation*}
f\in \mathcal{L}_{q,1,v}\Rightarrow T_{q,x}^{v}f\in \mathcal{L}_{q,1,v}.
\end{equation*}%
In [2] and as a first preamble of this theme the authors proved that
\begin{eqnarray*}
D_{-\frac{1}{2}}(q^{m},q^{r},q^{k}) &=&\frac{q^{2(r-m)(k-m)-m}}{(1-q)(q;q)_{I}}%
(q^{2(r-m)+1};q)_{\I}{_{1}\phi
_{1}}(0,q^{2(r-m)+1},q;q^{2(k-m)+1})\newline
\\
&=&\frac{1}{1-q}q^{-m}J_{2(r-m)}(q^{k-m},q),
\end{eqnarray*}%
which implies that the correspondent domain of $T_{q,x}^{-\frac{1}{2}}$ is
given by
\begin{equation*}
Q_{-\frac{1}{2}}=]0,q_{0}],
\end{equation*}%
where $q_{0}$ is the first zero of the $q$-hypergeometric function:
\begin{equation*}
q\mapsto {_{1}\phi _{1}(0,q,q,q)};
\end{equation*}%
a second statement is easily given by the Proposition 2 with $x=v=0$ which
gives
\begin{align*}
E_{0,0}(q^{m},q^{z},q^{k})& =c_{q,0}\int_{0}^{\infty }j_{0}(q^{m}\lambda
,q^{2})j_{0}(q^{z}\lambda ,q^{2})j_{0}(q^{k}\lambda ,q^{2})\lambda
d_{q}\lambda  \\
& =\frac{1}{c_{q,0}}D_{0}(q^{m},q^{z},q^{k}),
\end{align*}%
and then
\begin{align*}
D_{0}(q^{m},q^{z},q^{k})& =c_{q,0}E_{0,0}(q^{m},q^{z},q^{k}) \\
& =\frac{1}{(1-q)}q^{-2k}\left[ J_{z-k}(q^{m-k},q^{2})\right] ^{2},
\end{align*}%
Hence we have
\begin{equation*}
Q_{0}=]0,1[.
\end{equation*}%
Now we explicit the kernel in the production formula in terms of
$D_{v}$.

\begin{proposition}
For $n,m,k\in\mathbb{Z}$ and $-1<v$ we have
\begin{equation*}
E_{v,v}(q^{m},q^{n},q^{k})=(1-q)\sum_{i=0}^{\infty}\frac{(q^{-2v},q^{2})_{i}%
}{(q^{2},q^{2})_{i}}q^{i(2+2v)}D_{v}(q^{m},q^{i+n},q^{k}).
\end{equation*}
\end{proposition}

\begin{proof}
From the following formula (see[7], $\S $5) ``which can be proved in
a straightforward way by substitution of the defining series for the
$q$-Bessel functions on both sides, by interchanging summations, and
by evaluating the q-binomial series which occurs"
\begin{equation*}
J_{x-v}(\lambda ,q^{2})=\lambda ^{-v}\sum_{i=0}^{\infty }\frac{%
(q^{-2v},q^{2})_{i}}{(q^{2},q^{2})_{i}}q^{i(2+x)}J_{x}(\lambda q^{i},q^{2}),
\end{equation*}
where $\Re (x-v)>-1$ and $\Re (x)>-1$ we obtain
\begin{equation*}
(1-q)c_{q,x-v}j_{x-v}(\lambda ,q^{2})=(1-q)c_{q,x}\sum_{i=0}^{\infty }\frac{%
(q^{-2v},q^{2})_{i}}{(q^{2},q^{2})_{i}}q^{i(2+2x)}j_{x}(\lambda q^{i},q^{2}).
\end{equation*}%
Put $x=v$ and change $\lambda $ by $q^{n}\lambda $ we obtain
\begin{equation*}
j_{0}(q^{n}\lambda ,q^{2})=(1-q)c_{q,v}\sum_{i=0}^{\infty }\frac{%
(q^{-2v},q^{2})_{i}}{(q^{2},q^{2})_{i}}q^{i(2+2v)}j_{v}(\lambda
q^{i+n},q^{2}),
\end{equation*}%
which gives from Proposition 2
\begin{align*}
E_{v,v}(q^{m},q^{n},q^{k})& =c_{q,v}\int_{0}^{\infty }j_{v}(q^{m}\lambda
,q^{2})j_{0}(q^{n}\lambda ,q^{2})j_{v}(q^{k}\lambda ,q^{2})\lambda
^{2v+1}d_{q}\lambda \\
& =(1-q)\sum_{i=0}^{\infty }\frac{(q^{-2v},q^{2})_{i}}{(q^{2},q^{2})_{i}}%
q^{i(2+2v)} \\
& \times \left[ c_{q,v}^{2}\int_{0}^{\infty }j_{v}(q^{m}\lambda
,q^{2})j_{v}(q^{i+n}\lambda ,q^{2})j_{v}(q^{k}\lambda ,q^{2})\lambda
^{2v+1}d_{q}\lambda \right] \\
& =(1-q)\sum_{i=0}^{\infty }\frac{(q^{-2v},q^{2})_{i}}{(q^{2},q^{2})_{i}}%
q^{i(2+2v)}D_{v}(q^{m},q^{i+n},q^{k}).
\end{align*}%
We justify the exchange of the signs sum and integral by
$$\aligned
&\sum_{i=0}^{\infty }\frac{|(q^{-2v},q^{2})_{i}|}{(q^{2},q^{2})_{i}}%
q^{i(2+2v)}\int_{0}^{\infty }\left\vert j_{v}(q^{m}\lambda
,q^{2})j_{v}(q^{i+n}\lambda ,q^{2})j_{v}(q^{k}\lambda ,q^{2})\lambda
^{2v+1}\right\vert d_{q}\lambda\\
&\leq\Vert j_{v}(.,q^{2})\Vert _{q,\infty }^{2}\Vert
j_{v}(.,q^{2})\Vert
_{q,1,v}q^{-2(v+1)m}\sum_{i=0}^{\infty }\frac{|(q^{-2v},q^{2})_{i}|}{%
(q^{2},q^{2})_{i}}q^{i(2+2v)}<\infty .
\endaligned$$
So we obtain the result.
\end{proof}

\begin{proposition}
Let $-1<v<0$ then
\begin{equation*}
Q_v\varsubsetneq ]0,1[.
\end{equation*}
\end{proposition}

\begin{proof}
For $m,n,k\in \mathbb{Z}$ and $-1<v<0$ we have%
\begin{eqnarray*}
E_{v,v}(q^{m},q^{n},q^{k}) &=&\frac{1}{(1-q)c_{q,v}}%
q^{-k(v+2)-mv}J_{n-k}(q^{m-k},q^{2})J_{v+n-k}(q^{m-k},q^{2})\newline
\\
&=&(1-q)\sum_{i=0}^{\infty }\frac{(q^{-2v},q^{2})_{i}}{(q^{2},q^{2})_{i}}%
q^{i(2+2v)}D_{v}(q^{m},q^{i+n},q^{k}).
\end{eqnarray*}%
In particular for $m=n=k=0$
\begin{equation*}
E_{v,v}(1,1,1)=\frac{1}{(1-q)c_{q,v}}J_{0}(1,q^{2})J_{v}(1,q^{2}).
\end{equation*}%
We introduce the following function
\begin{equation*}
\phi _{v}:q\mapsto (q^{2},q^{2})_{\infty }^{2}J_{v}(1,q^{2}).
\end{equation*}%
Let $q_{1}\simeq 0.658$ the first zero of $\phi _{0}$ and consider the graph
of the function
\begin{equation*}
v\mapsto \phi _{v}(q_{1}),\quad v\in \lbrack 0,1].
\end{equation*}
$$
\scalebox{0.5}{\epsfbox{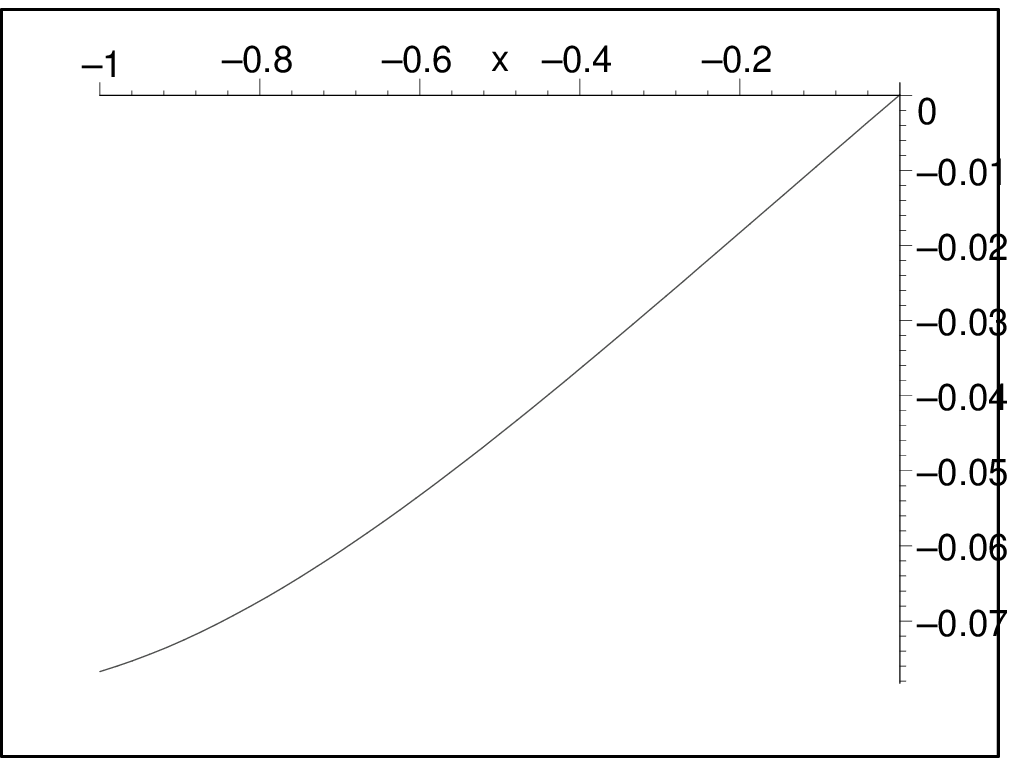}}
$$
We conclude that for $-1<v<0$ we have $\phi _{v}(q_{1})<0$. Then
there exist a very small $\varepsilon >0$ such that $\phi
_{v}(q_{1}-\varepsilon )<0$ and $\phi _{0}(q_{1}-\varepsilon )>0$.
Hence, this function
\begin{equation*}
q\mapsto J_{0}(1,q^{2})J_{v}(1,q^{2})<0,
\end{equation*}%
takes some negative values in the interval $]0,1[$, then for some entire $%
i\in \mathbb{N}$
\begin{equation*}
\frac{(q^{-2v},q^{2})_{i}}{(q^{2},q^{2})_{i}}q^{i(2+2v)}D_{v}(1,q^{i},1)<0.
\end{equation*}%
As
\begin{equation*}
\frac{(q^{-2v},q^{2})_{i}}{(q^{2},q^{2})_{i}}q^{i(2+2v)}>0,
\end{equation*}%
then $D_{v}(1,q^{i},1)<0$ for some entire $i\in \mathbb{N}$. Which prove
that
\begin{equation*}
Q_{v}\varsubsetneq ]0,1[.
\end{equation*}%
This finish the proof.
\end{proof}

\begin{lemma}
For $x\in\mathbb{R}_q^+$ and $\Re(v+t)>-1$ we have
\begin{equation*}
(1-q)c_{q,v}\sum_{n\in\mathbb{Z}}q^{(1+v+t)n}j_v(q^nx,q^2)=\frac{%
(q^{1+v-t};q^2)_\I}{(q^{1+v+t};q^2)_\I}x^{-(1+v+t)}.
\end{equation*}
\end{lemma}

\begin{proof}
In [6] the following result was proved
\begin{equation*}
\sum_{n\in \mathbb{Z}}q^{(1+t)n}J_{v}(q^{n},q^{2})=\frac{(q^{1+v-t};q^{2})_{%
\infty }}{(q^{1+v+t};q^{2})_{\infty }},\quad \Re (v+t)>-1.
\end{equation*}%
Then we get
\begin{equation*}
(1-q)c_{q,v}\sum_{n\in \mathbb{Z}}q^{(1+v+t)n}j_{v}(q^{n},q^{2})=\frac{%
(q^{1+v-t};q^{2})_{\infty }}{(q^{1+v+t};q^{2})_{\infty }}.
\end{equation*}%
Hence for $x=q^{k}\in \mathbb{R}_{q}^{+}$ we have
\begin{equation*}
(1-q)c_{q,v}\sum_{n\in \mathbb{Z}}q^{(1+v+t)(n+k)}j_{v}(q^{n+k},q^{2})=\frac{%
(q^{1+v-t};q^{2})_{\infty }}{(q^{1+v+t};q^{2})_{\infty }}.
\end{equation*}%
This finish the proof.
\end{proof}

\begin{theorem}
Let $v\geq x>-1$ then
\begin{equation*}
Q_{x}\subset Q_{v}.
\end{equation*}%
As consequence

\begin{itemize}
\item If $0\leq v$ then $Q_{v}=]0,1[.$

\item If $-\frac{1}{2}\leq v<0$ then $]0,q_{0}]\subset $ $Q_{v}\subsetneq
]0,1[.$

\item If $-1<v\leq -\frac{1}{2}$ then $Q_{v}\subset ]0,q_{0}].$
\end{itemize}
\end{theorem}

\begin{proof}
Let $v> -1$ and $0<\mu <1$. We will prove that
\begin{equation*}
Q_{v}\subset Q_{v+\mu },
\end{equation*}%
which suffices to show that $Q_{x}\subset Q_{v}$ if $x\leq v$. In
the following we tack $q\in Q_{v}.$ We have
\begin{equation*}
c_{q,v+\mu }j_{v+\mu }(t,q^{2})=c_{q,v}\sum_{i=0}^{\infty }\frac{(q^{2\mu
},q^{2})_{i}}{(q^{2},q^{2})_{i}}q^{2(v+1)i}j_{v}(tq^{i},q^{2}),\newline
\end{equation*}%
and then%
$$\aligned
&c_{q,v}c_{q,v+\mu }^{2}j_{v}(tx,q^{2})j_{v+\mu }(ty,q^{2})j_{v}(tz,q^{2})\\
&=c_{q,v}^{3}\sum_{i,j,s=0}^{\infty }\frac{(q^{2\mu },q^{2})_{i}}{%
(q^{2},q^{2})_{i}}\frac{(q^{2\mu },q^{2})_{j}}{(q^{2},q^{2})_{j}}%
q^{2(1+v)(i+j)}\\
&\times j_{v}(tx,q^{2})j_{v}(tq^{i}y,q^{2})j_{v}(tq^{j}z,q^{2}),
\endaligned$$
which implies
$$\aligned
&c_{q,v}T_{v,v+\mu ,v}(x,y,z)\\
&=c_{q,v}\sum_{i,j=0}^{\infty }\frac{(q^{2\mu },q^{2})_{i}}{%
(q^{2},q^{2})_{i}}\frac{(q^{2\mu },q^{2})_{j}}{(q^{2},q^{2})_{j}}%
q^{2(1+v)(i+j)}D_{v}(x,q^{i}y,q^{j}z)\geq 0,
\endaligned$$
where
\begin{equation*}
T_{v,w,\alpha }(x,y,z)=c_{q,w}^{2}\int_{0}^{\infty
}j_{v}(tx,q^{2})j_{w}(ty,q^{2})j_{w}(tz,q^{2})t^{2\alpha +1}d_{q}t.
\end{equation*}%
Note that
\begin{equation*}
T_{v,v,v}(x,y,z)=D_{v}(x,y,z).
\end{equation*}%
The exchange of the signs sum and integral is valid since:
\begin{align*}
& \sum_{i,j=0}^{\infty }q^{2(1+v)(i+j)}\int_{0}^{\infty }\left\vert
j_{v}(q^{m}t,q^{2})\right\vert \left\vert j_{v}(q^{n+i}t,q^{2})\right\vert
\left\vert j_{v}(q^{k+j}t,q^{2})\right\vert t^{2v+1}d_{q}t \\
& =\sum_{i,j=0}^{\infty }q^{2(1+v)(i+j)}\int_{0}^{1}\left\vert
j_{v}(q^{m}t,q^{2})\right\vert \left\vert j_{v}(q^{n+i}t,q^{2})\right\vert
\left\vert j_{v}(q^{k+j}t,q^{2})\right\vert t^{2v+1}d_{q}t \\
& +\sum_{i,j=0}^{\infty }q^{2(1+v)(i+j)}\int_{1}^{\infty }\left\vert
j_{v}(q^{m}t,q^{2})\right\vert \left\vert j_{v}(q^{n+i}t,q^{2})\right\vert
\left\vert j_{v}(q^{k+j}t,q^{2})\right\vert t^{2v+1}d_{q}t.
\end{align*}%
Note that
\begin{align*}
& \sum_{i,j=0}^{\infty }q^{2(1+v)(i+j)}\int_{0}^{1}\left\vert
j_{v}(q^{m}t,q^{2})\right\vert \left\vert j_{v}(q^{n+i}t,q^{2})\right\vert
\left\vert j_{v}(q^{k+j}t,q^{2})\right\vert t^{2v+1}d_{q}t \\
& \leq \left\Vert j_{v}(.,q^{2})\right\Vert _{q,\infty }^{3}\left[
\sum_{i,j=0}^{\infty }q^{2(1+v)(i+j)}\right]
\int_{0}^{1}t^{2v+1}d_{q}t<\infty ,
\end{align*}%
and
\begin{align*}
& \sum_{i,j=0}^{\infty }q^{2(1+v)(i+j)}\int_{1}^{\infty }\left\vert
j_{v}(q^{m}t,q^{2})\right\vert \left\vert j_{v}(q^{n+i}t,q^{2})\right\vert
\left\vert j_{v}(q^{k+j}t,q^{2})\right\vert t^{2v+1}d_{q}t \\
& =(1-q)\sum_{i,j,r^{\prime }=0}^{\infty }q^{2(1+v)(i+j-r)}\left\vert
j_{v}(q^{m-r},q^{2})\right\vert \left\vert j_{v}(q^{n+i-r},q^{2})\right\vert
\left\vert j_{v}(q^{k+j-r},q^{2})\right\vert  \\
& =\sum_{j,r^{\prime }=0}^{\infty }q^{2(1+v)j}\left\vert
j_{v}(q^{m-r},q^{2})\right\vert \left\vert
j_{v}(q^{k+j-r},q^{2})\right\vert \times \left[
(1-q)\sum_{i=0}^{\infty }q^{2(1+v)(i-r)}\left\vert
j_{v}(q^{n+i-r},q^{2})\right\vert \right],
\end{align*}%
where $r^{\prime }=r-1.$ We write
\begin{align*}
\sum_{i=0}^{\infty }q^{2(1+v)(i-r)}\left\vert
j_{v}(q^{n+i-r},q^{2})\right\vert & =
q^{-2(1+v)n}(1-q)\sum_{i=0}^{\infty }q^{2(n+i-r)}\left\vert
j_{v}(q^{n+i-r},q^{2})\right\vert   \\
& =q^{-2(1+v)n}(1-q)\sum_{i=m-r}^{\infty }q^{2(1+v)i}\left\vert
j_{v}(q^{i},q^{2})\right\vert   \\
& <q^{-2(1+v)n}\left\Vert j_{v}(.,q^{2})\right\Vert _{q,1,v}.
\end{align*}%
Then
\begin{align*}
& \sum_{j,r^{\prime }=0}^{\infty }q^{2(1+v)j}\left\vert
j_{v}(q^{m-r},q^{2})\right\vert \left\vert j_{v}(q^{k+j-r},q^{2})\right\vert
\\
& =\sum_{s,r^{\prime }=0}^{\infty }\left\vert
j_{v}(q^{m-r},q^{2})\right\vert \left[ \sum_{j=0}^{\infty
}q^{2(1+v)j}\left\vert j_{v}(q^{k+j-r},q^{2})\right\vert \right]  \\
& =\sum_{r^{\prime }=0}^{\infty }\left\vert
j_{v}(q^{m-r},q^{2})\right\vert \left[
q^{2(1+v)(r-k)}\sum_{j=0}^{\infty }q^{2(k+j-r)}\left\vert
j_{0}(q^{k+j-r},q^{2})\right\vert \right]  \\
& =\sum_{s,r^{\prime }=0}^{\infty }q^{2s}\left\vert
j_{v}(q^{m-r},q^{2})\right\vert \left[
q^{2(v+1)(r-k)}\sum_{j=k-r}^{\infty
}q^{2(1+v)j}\left\vert j_{0}(q^{j},q^{2})\right\vert \right]  \\
& \leq \frac{q^{-2(v+1)k}}{1-q}\left\Vert j_{v}(.,q^{2})\right\Vert
_{q,1,v}\left\Vert j_{v}(.,q^{2})\right\Vert _{q,\infty
}\sum_{r^{\prime }=0}^{\infty }q^{2(v+1)r}<\infty .
\end{align*}

Let $0<\alpha <\mu $. We introduce the function $A_{\alpha ,\mu ,v}(x)$ as
follows:
\begin{equation*}
A_{\alpha ,\mu ,v}(x)=c_{q,v}\int_{0}^{\infty }t^{2\alpha }j_{v+\mu
}(t,q^{2})j_{v}(xt,q^{2})t^{2v+1}d_{q}t.
\end{equation*}%
From the inversion formula in Theorem 1 ($t\mapsto t^{2\alpha }j_{v+\mu
}(t,q^{2})\in \mathcal{L}_{q,1,v}$) we get
\begin{equation*}
c_{q,v}\int_{0}^{\infty }A_{\alpha ,\mu
,v}(x)j_{v}(xt,q^{2})x^{2v+1}d_{q}x=t^{2\alpha }j_{v+\mu }(t,q^{2}).
\end{equation*}%
Let $x\leq 1$. From Lemma 1 we obtain
\begin{align*}
& A_{\alpha ,\mu ,v}(x) \\
& =c_{q,v}(1-q)\sum_{n\in \mathbb{Z}}q^{2(\alpha
+v+1)n}j_{v}(xq^{n},q^{2})j_{v+\mu }(q^{n},q^{2}) \\
& =c_{q,v}(1-q)\sum_{n\in \mathbb{Z}}\left[ \sum_{i=0}^{\infty }(-1)^{i}%
\frac{q^{i(i+1)}}{(q^{2+2v},q^{2})_{i}(q^{2},q^{2})_{i}}\right] q^{2(\alpha
+v+1+i)n}j_{v+\mu }(q^{n},q^{2}) \\
& =c_{q,v}(1-q)\sum_{i=0}^{\infty }(-1)^{i}\frac{q^{i(i+1)}}{%
(q^{2+2v},q^{2})_{i}(q^{2},q^{2})_{i}}x^{2i}\left[ \sum_{n\in \mathbb{Z}%
}q^{2(\alpha +v+1+i)n}j_{v+\mu }(q^{n},q^{2})\right]  \\
& =\frac{c_{q,v}}{c_{q,v+\mu }}\sum_{i=0}^{\infty }(-1)^{i}\frac{q^{i(i+1)}}{%
(q^{2+2w},q^{2})_{i}(q^{2},q^{2})_{i}}x^{2i}\frac{(q^{2(1+v+\mu )-2(\alpha
+v+1+i)},q^{2})_{\infty }}{(q^{2(\alpha +v+1+i)},q^{2})_{\infty }} \\
& =\frac{c_{q,v}}{c_{q,v+\mu }}\sum_{i=0}^{\infty }(-1)^{i}\frac{q^{i(i+1)}}{%
(q^{2+2v},q^{2})_{i}(q^{2},q^{2})_{i}}x^{2i}\frac{(q^{2(\mu -\alpha
)}q^{-2i},q^{2})_{\infty }}{(q^{2(\alpha +v+1+i)},q^{2})_{\infty }}\geq 0.
\end{align*}%
Note that
\begin{equation*}
0<\mu -\alpha <1.
\end{equation*}%
The exchange of the signs sum is valid since. Indeed, let $q^{k}=x$ then we
have
\begin{align*}
& \sum_{i=0}^{\infty }q^{i(i+1)}q^{2ik}\left[ \sum_{n\in \mathbb{Z}%
}q^{2(\alpha +v+1+i)n}\left\vert j_{v+\mu }(q^{n},q^{2})\right\vert \right]
\\
& =\sum_{i=0}^{\infty }q^{i(i+1)}q^{2ik}\left[ \sum_{n=0}^{\infty
}q^{2(\alpha +v+1+i)n}\left\vert j_{v+\mu }(q^{n},q^{2})\right\vert \right]
\\
& +\sum_{i=0}^{\infty }q^{i(i+1)}q^{2ik}\left[ \sum_{n=1}^{\infty
}q^{-2(\alpha +v+1+i)n}\left\vert j_{v+\mu }(q^{-n},q^{2})\right\vert \right]
\\
& \leq \sum_{i=0}^{\infty }q^{i(i+1)}q^{2ik}\left[ \sum_{n=0}^{\infty
}q^{2(\alpha +v+1)n}\right]  \\
& +\sum_{i=0}^{\infty }q^{i(i+1)}q^{2ik}\left[ \sum_{n=1}^{\infty
}q^{-2(\alpha +v+1+i)n+n^{2}+(2v+2\mu +1)n}\right] .
\end{align*}%
One has to observe that the first sum exist. The second sum also:
\begin{align*}
& \sum_{i=0}^{\infty }q^{i(i+1)}q^{2ik}\left[ \sum_{n=1}^{\infty
}q^{-2(\alpha +v+1+i)n+n^{2}+(2v+2\mu +1)n}\right]  \\
& =\sum_{n=1}^{\infty }q^{-2(\alpha +v+1+i)n+n^{2}+(2v+2\mu +1)n}\left[
\sum_{i=0}^{\infty }q^{i(i+1)}q^{2ik}\right]  \\
& =\sum_{n=1}^{\infty }q^{-2(\alpha +v+1)n+(2v+2\mu +1)n+(2k+1)n}\left[
\sum_{i=0}^{\infty }q^{(i-n)^{2}}q^{(2k+1)(i-n)}\right]  \\
& =\sum_{n=1}^{\infty }q^{2(\mu -\alpha )n+(2k+1)n}\left[ \sum_{i=-n}^{%
\infty }q^{i^{2}}q^{(2k+1)i}\right] \leq \sum_{n=1}^{\infty }q^{2(\mu
-\alpha )n+(2k+1)n}\left[ \sum_{i=-\infty }^{\infty }q^{i^{2}}q^{(2k+1)i}%
\right] <\infty .
\end{align*}%
If $x>1$ then we obtain%
\begin{align*}
& A_{\alpha ,\mu ,v}(x) \\
& =c_{q,v}(1-q)\sum_{n\in \mathbb{Z}}q^{2(\alpha
+v+1)n}j_{v}(xq^{n},q^{2})j_{v+\mu }(q^{n},q^{2}) \\
& =c_{q,v}(1-q)\sum_{n\in \mathbb{Z}}\left[ \sum_{i=0}^{\infty }(-1)^{i}%
\frac{q^{i(i+1)}}{(q^{2+2v},q^{2})_{i}(q^{2},q^{2})_{i}}\right] q^{2(\alpha
+v+1+i)n}j_{v}(q^{n}x,q^{2}) \\
& =c_{q,v}(1-q)\sum_{i=0}^{\infty }(-1)^{i}\frac{q^{i(i+1)}}{(q^{2+2v+2\mu
},q^{2})_{i}(q^{2},q^{2})_{i}}\left[ \sum_{n\in \mathbb{Z}}q^{2(\alpha
+v+1+i)n}j_{v}(q^{n}x,q^{2})\right]  \\
& =x^{-2(\alpha +v+1)}\sum_{i=0}^{\infty }(-1)^{i}\frac{q^{i(i+1)}}{%
(q^{2+2v+2\mu },q^{2})_{i}(q^{2},q^{2})_{i}}x^{-2i}\frac{(q^{2(1+v)-2(\alpha
+v+1+i)},q^{2})_{\infty }}{(q^{2(\alpha +v+1+i)},q^{2})_{\infty }} \\
& =x^{-2(\alpha +v+1)}\sum_{i=0}^{\infty }(-1)^{i}\frac{q^{i(i+1)}}{%
(q^{2+2v+2\mu },q^{2})_{i}(q^{2},q^{2})_{i}}x^{-2i}\frac{(q^{-2\alpha
}q^{-2i},q^{2})_{\infty }}{(q^{2(\alpha +v+1+i)},q^{2})_{\infty }}<0.
\end{align*}%
Now, we write
\begin{align*}
& c_{q,v}\int_{0}^{\infty }A_{\alpha ,\mu ,v}(x)T_{v,v+\mu
,v}(x,y,z)x^{2v+1}d_{q}x \\
& =c_{q,v}\int_{0}^{\infty }A_{\alpha ,\mu ,v}(x)\left[ c_{q,v}^{2}%
\int_{0}^{\infty }j_{w}(xt,q^{2})j_{v}(yt,q^{2})j_{v}(zt,q^{2})t^{2v+1}d_{q}t%
\right] x^{2v+1}d_{q}x \\
& =c_{q,v+\mu }^{2}\int_{0}^{\infty }\left[ c_{q,v}\int_{0}^{\infty
}A_{\alpha ,\mu ,v}(x)j_{v}(xt,q^{2})x^{2v+1}d_{q}x\right] j_{v+\mu
}(yt,q^{2})j_{v+\mu }(zt,q^{2})t^{2v+1}d_{q}t \\
& =c_{q,v+\mu }^{2}\int_{0}^{\infty }j_{v+\mu }(t,q^{2})j_{v+\mu
}(yt,q^{2})j_{v+\mu }(zt,q^{2})t^{2(v+\alpha )+1}d_{q}t=T_{v+\mu ,v+\mu
,v+\alpha }(1,y,z).
\end{align*}%
To justify the exchange of the signs integrals we write
$$\aligned
&\int_{0}^{\infty }\left[ \int_{0}^{\infty }|A_{\alpha ,\mu
,v}(x)||j_{v}(xt,q^{2})|x^{2v+1}d_{q}x\right]
|j_{v}(yt,q^{2})||j_{v}(zt,q^{2})|td_{q}t\\
&\leq\Big[\Vert j_{v}(.,q^{2})\Vert_{q,\infty }^{2}\Vert
j_{v}(.,q^{2})\Vert _{q,1,v}\Vert A_{\alpha ,\mu ,v}\Vert
_{q,1,v}\Big]\frac{1}{z}.
\endaligned$$
Not that $x\mapsto A_{\alpha ,\mu ,v}(x)$ is continued at $0$ and
\begin{equation*}
|A_{\alpha ,\mu ,v}(x)|\leq \Big[\Vert j_{v+\mu }(.,q^{2})\Vert _{q,\infty
}\Vert j_{v}(.,q^{2})\Vert _{q,1,\alpha +v}\Big]x^{-2(\alpha +v+1)},\quad
\text{as}\quad x\rightarrow \infty .
\end{equation*}%
On the other hand%
\begin{equation*}
\int_{0}^{\infty }A_{\alpha ,\mu ,v}(x)x^{2v+1}d_{q}x=0,
\end{equation*}%
then we obtain for all $\delta >0$%
\begin{equation*}
T_{v+\mu ,v+\mu ,v+\alpha }(1,y,z)=\int_{0}^{\infty }A_{\alpha ,\mu
,v}\left( x\right) \Big[T_{v,v+\mu ,v}(x,y,z)-\delta \Big]x^{2v+1}d_{q}x%
\newline
.
\end{equation*}%
In the following we assume that $0<y,z\leq 1.$ Let
\begin{equation*}
\delta _{0}=\inf_{0\leq y,z\leq 1}T_{v,v+\mu ,v}(1,y,z)
\end{equation*}%
Not that $\delta _{0}$ exist and strictly positive, indeed the following
function%
\begin{equation*}
(y,z)\mapsto T_{v,v+\mu ,v}(1,y,z)
\end{equation*}%
is continuous on the compact $[0,1]\times \lbrack 0,1].$ Hence,
there exist
\begin{equation*}
(y_{0},z_{0})\in \lbrack 0,1]\times \lbrack 0,1]
\end{equation*}%
such that%
\begin{equation*}
\delta _{0}=T_{v,v+\mu ,v}(1,y_{0},z_{0})=\inf_{0\leq y,z\leq
1}T_{v,v+\mu ,v}(1,y,z)\geq 0.
\end{equation*}%
If we assume that $T_{v,v+\mu ,v}(1,y_{0},z_{0})=0$ $\ $then%
\begin{equation*}
\sum_{i,j=0}^{\infty }\frac{(q^{2\mu },q^{2})_{i}}{(q^{2},q^{2})_{i}}\frac{%
(q^{2\mu },q^{2})_{j}}{(q^{2},q^{2})_{j}}%
q^{2(1+v)(i+j)}D_{v}(1,q^{i}y_{0},q^{j}z_{0})=0,
\end{equation*}%
which implies that%
\begin{equation*}
c_{q,v}\mathcal{F}_{q,v}\Big[t\mapsto j_{v}(t,q^{2})j_{v}(z_{0}t,q^{2})\Big]%
(q^{i}y_{0})=D_{v}(1,q^{i}y_{0},z_{0})=0,\text{ }\forall i\in \mathbb{N}.
\end{equation*}%
From Proposition 1 and the fact that there exist $\sigma_0>0$ such
that
\begin{equation*}
\left\vert j_{v}(z,q^{2})\right\vert \leq \sigma _{0}\text{ }e^{\left\vert
z\right\vert },\text{ \ \ \ \ \ \ \ \ }\forall z\in \mathbb{C}
\end{equation*}%
we see that this function
\begin{equation*}
z\mapsto \mathcal{F}_{q,v}\Big[t\mapsto j_{v}(t,q^{2})j_{v}(z_{0}t,q^{2})%
\Big](z)
\end{equation*}%
is analytic. Then for all $x\in \mathbb{R}_{q}$%
\begin{equation*}
\mathcal{F}_{q,v}\Big[t\mapsto j_{v}(t,q^{2})j_{v}(z_{0}t,q^{2})\Big]%
(x)=0\Rightarrow j_{v}(t,q^{2})j_{v}(z_{0}t,q^{2})=0,\text{ \ }\forall t\in
\mathbb{R}_{q},
\end{equation*}%
but this is absurd. Then $\delta _{0}>0.$ Now we have
\begin{equation*}
\delta _{0}\leq T_{v,v+\mu ,v}(x,y,z),\text{ \ \ \ }\forall 0<x\leq 1.
\end{equation*}%
If
\begin{equation*}
\delta _{0}>T_{v,v+\mu ,v}(x,y,z),\text{ \ \ \ }\forall x>1.
\end{equation*}%
then
\begin{eqnarray*}
T_{v+\mu ,v+\mu ,v+\alpha }(1,y,z) &=&\int_{0}^{1}A_{\alpha ,\mu ,v}\left(
x\right) \Big[T_{v,v+\mu ,v}(x,y,z)-\delta _{0}\Big]x^{2v+1}d_{q}x\newline
\\
&&+\int_{1}^{\infty }A_{\alpha ,\mu ,v}\left( x\right) \Big[T_{v,v+\mu
,v}(x,y,z)-\delta _{0}\Big]x^{2v+1}d_{q}x\geq 0.
\end{eqnarray*}%
Otherwise, there exist $s>0$ such that
\begin{equation*}
\delta _{0}>T_{v,v+\mu ,v}(x,y,z),\text{ \ \ \ }\forall x>s.
\end{equation*}%
because%
\begin{equation*}
\left\vert T_{v,v+\mu ,v}(x,y,z)\right\vert <c_{q,v+\mu }^{2}\left\Vert
j_{v+\mu }(.q^{2})\right\Vert _{q,\infty }^{2}\times \left\Vert
j_{v}(.q^{2})\right\Vert _{q,v,1}x^{-2(v+1)}\rightarrow 0,\text{ \ \ \ as }%
x\rightarrow \infty .
\end{equation*}%
For $\delta \geq \delta _{0}$ we obtain
$$\aligned
T_{v+\mu ,v+\mu ,v+\alpha }(1,y,z)&=\int_{0}^{1}A_{\alpha ,\mu
,v}\left( x\right) \Big[T_{v,v+\mu ,v}(x,y,z)-\delta
\Big]x^{2v+1}d_{q}x\\
&+\int_{1}^{s}A_{\alpha ,\mu ,v}\left( x\right) \Big[T_{v,v+\mu
,v}(x,y,z)-\delta \Big]x^{2v+1}d_{q}x\\
&+\int_{s}^{\infty }A_{\alpha ,\mu ,v}\left( x\right)
\Big[T_{v,v+\mu
,v}(x,y,z)-\delta \Big]x^{2v+1}d_{q}x \\
&=I_{1}(\delta )+I_{2}(\delta )+I_{3}(\delta ).
\endaligned$$
Note that%
\begin{equation*}
\left\{
\begin{tabular}{l}
$\delta \mapsto I_{1}(\delta )$ is a decreesing function tends towards $\
-\infty $ and $I_{1}(\delta _{0})>0$ \\
$\delta \mapsto I_{2}(\delta )$ is an increasing function tends towards $%
+\infty $ and $I_{2}(\delta _{0})<0$ \\
$\delta \mapsto I_{3}(\delta )$ is an increasing function tends towards $%
+\infty $ \ and $I_{3}(\delta _{0})>0$%
\end{tabular}%
\right.
\end{equation*}%
then there exist $\delta >\delta _{0\text{ }}$ such that
\begin{equation*}
I_{1}(\delta )+I_{2}(\delta )=0
\end{equation*}%
which implies%
\begin{equation*}
T_{v+\mu ,v+\mu ,v+\alpha }(1,y,z)=I_{3}(\delta )>0.
\end{equation*}%
In the end%
$$\aligned
&\lim_{\alpha \rightarrow \mu }\left[ T_{v+\mu ,v+\mu ,v+\alpha
}(1,y,z)\right]\\
&=\lim_{\alpha \rightarrow \mu }\left[ c_{q,v+\mu
}^{2}\int_{0}^{\infty }j_{v+\mu }(t,q^{2})j_{v+\mu
}(yt,q^{2})j_{v+\mu }(zt,q^{2})t^{2(v+\alpha
)+1}d_{q}t\right]\\
&=c_{q,v+\mu }^{2}\int_{0}^{\infty }j_{v+\mu }(t,q^{2})j_{v+\mu
}(yt,q^{2})j_{v+\mu }(zt,q^{2})\big(\lim_{\alpha \rightarrow \mu
}t^{2(v+\alpha )+1}\Big)d_{q}t\\
&=T_{v+\mu ,v+\mu ,v+\alpha }(1,y,z)=D_{v+\mu }(1,y,z)\geq 0.
\endaligned$$
It is not hard to justify the exchange of the signs integral and
limit, in
fact%
$$\aligned
&\left\vert j_{v+\mu }(t,q^{2})j_{v+\mu }(yt,q^{2})j_{v+\mu
}(zt,q^{2})t^{2(v+\alpha )+1}\right\vert\\
&\leq\left\Vert t\mapsto j_{v+\mu }(yt,q^{2})j_{v+\mu
}(zt,q^{2})t^{2(v+\alpha )+1}\right\Vert _{q,\infty }\times
\left\vert j_{v+\mu }(t,q^{2})\right\vert .
\endaligned$$
From the following identity
\begin{equation*}
D_{v+\mu }(x,y,z)=x^{-2(v+\mu +1)}D_{v+\mu }\left( 1,\frac{y}{x},\frac{z}{x}%
\right) ,
\end{equation*}%
we deduce that%
\begin{equation*}
D_{v+\mu }(x,y,z)\geq 0,\text{ \ \ \ \ \ \ \ }\forall x,y,z\in \mathbb{R}%
_{q}^{+}.
\end{equation*}%
The fact that
\begin{equation*}
Q_{-\frac{1}{2}}=]0,q_{0}]\text{ \ and }Q_{0}=]0,1[
\end{equation*}%
leads to the result.
\end{proof}

\end{document}